\input amstex
\documentstyle{amsppt}
\magnification=1200
\pageheight{7.5in}
\expandafter\redefine\csname logo\string@\endcsname{}
\NoBlackBoxes

\topmatter
\title Correction to `$K$-theory of virtually poly-surface groups' 
\endtitle
\author S. K. Roushon
\endauthor
\address School of Mathematics, Tata Institute, Homi Bhabha Road, Mumbai
400 005, India. 
\endaddress
\email roushon\@math.tifr.res.in \newline
\indent {\it URL:}\ http://www.math.tifr.res.in/\~\
roushon/paper.html \endemail
\date August 18, 2004 \enddate
\abstract In this note we point out an error in the above
paper and refer to some papers where this error is corrected and a
more general theorem is proved.
\endabstract
\keywords Strongly poly-surface groups, Whitehead group, fibered
isomorphism conjecture, negative $K$-groups
\endkeywords
\subjclass  Primary: 19B28, 19A31, 20F99, 19D35. Secondary: 19J10
\endsubjclass
\endtopmatter
\document
\baselineskip 14pt
In this note `FIC' stands for the Fibered Isomorphism
Conjecture of Farrell and Jones corresponding to the pseudoisotopy
functor (see \cite{FJ}).

In the proof of the main lemma of \cite{R1} we found some
filtration of the surface $\tilde F$ which is preserved by the
diffeomorphism $f$ and
used this filtration to find a filtration of the mapping torus $M_f$ of
$f$ by compact submanifolds with incompressible tori boundary. Recall
that $\tilde F$ was the covering of the surface $F$
corresponding to the commutator subgroup of $\pi_1(F)$ and $f:\tilde F\to 
\tilde F$ was a lift of a diffeomorphism $g:F\to F$. Also recall that the
main lemma of \cite{R1} says that the FIC is true for $\pi_1(M)$ where
$M$ is the mapping torus of a diffeomorphism of $F$. The proof of the
existence of the above filtration of $M_f$, we sketched in \cite{R1} is
incorrect. In the proof of the main lemma of \cite{R2} we show that some
regular finite sheeted cover of $M_f$ admits a filtration of the required
type provided $g$ satisfies certain conditions. For general $g$ we
prove the main lemma of \cite{R1} in \cite{R3} assuming that the FIC
is
true for $B$-groups. By definition a $B$-group contains a finite index 
subgroup isomorphic to the fundamental group of a compact irreducible
$3$-manifold with nonempty incompressible boundary so that each boundary
component is a surface of genus $\geq 2$. Finally we refer to [\cite{R4},
corollary 1.1] 	where we prove that the FIC is true for $B$-groups. This
completes the proof of the main lemma of \cite{R1}. 

Recall that the
strongly poly-surface groups ([\cite{R1}, definition]) were defined in
\cite{R1} and the FIC was proved for any virtually strongly
poly-surface group ([\cite{R1}, main theorem]). We
generalized this notion and defined weak strongly
poly-surface groups (see [\cite{R4}, definition 1.1]) and in [\cite{R4},
theorem 1.2] the FIC is proved for any virtually weak strongly
poly-surface group.  

Also we should point out that in this situation the proof of
proposition 2.3 of \cite{R1} needs a slightly elaborate argument.
This proof is now contained in the proof of the main theorem of
\cite{R3} and stated in corollary 3.5 of \cite{R3}. Recall that
proposition 2.3 says that the FIC is true for the fundamental group of a
$3$-manifold which has a finite sheeted cover fibering over the circle.

\enddocument